# Two definitions of fractional derivatives of powers functions

**Raoelina Andriambololona, Rakotoson Hanitriarivo, Tokiniaina Ranaivoson, Roland Raboanary**

Theoretical Physics Dept., Antananarivo, Madagascar
Institut National des Sciences et Techniques Nucléaires (INSTN-Madagascar), Antananarivo, Madagascar

**Email address:**
raoelinasp@yahoo.fr (R. Andriambololona), infotsara@gmail.com (R. Hanitriarivo), tokhiniaina@gmail.com (T. Ranaivoson)

**Abstract:** We consider the set of powers functions defined on $\mathbb{R}_+$ and their linear combinations. After recalling some properties of the gamma function, we give two general definitions of derivatives of positive and negative integers, positive and negative fractional orders. Properties of linearity and commutativity are studied and the notions of semi-equality, semi-linearity and semi-commutativity are introduced. Our approach gives a unified definition of the common derivatives and integrals and their generalization.

**Keywords:** Gamma Function; Fractional Derivatives; Fractional Integrals; Power Functions

## 1. Introduction

The concept of fractional derivative can be considered as rising from the idea to make sense to the expression

$$f^{(s)}(x) = \frac{d^s f(x)}{dx^s} \qquad (1.1)$$

for a function $f$ of a variable $x$ and for a fractional number $s$. Historically, the origin of the concept of fractional derivative was attributed to Leibnitz and l'Hospital (1695) but significant developments on the subject were given later on by Riemann, Liouville, Grunwald-Letnikov and others. Several approaches are considered by various authors in the definitions of this concept [1], [2], [3], [4], [5] [6], [7]. In our case, the study is based on the definition of fractional derivative for power functions defined on the set $\mathbb{R}_+$ of positive number and their linear combination. We utilize the properties of the gamma functions [8],[9] [10], [11].

One of the major results that we obtain is the unification of the definition for both positive and negative order of the fractional derivative.

## 2. Recall on Some Properties of the Gamma Function

The gamma function, noted $\Gamma$, is defined as the integral

$$\Gamma(z) = \int_0^{+\infty} t^{z-1} e^{-t} \, dt \qquad (2.1)$$

for $z \in \mathbb{N}$ and $t \in ]0 + \infty[$. The following properties are easily deduced

$$\begin{cases} \Gamma(1) = 1 \\ \Gamma(z) = \dfrac{\Gamma(z+1)}{z} \end{cases} \qquad (2.2)$$

$$\Gamma(z) = (z-1)! \qquad (2.3)$$

The definition of $\Gamma$ may be generalized starting from the relations (2.1) and (2.2) to $z$ negative integer, fractional, real and complex. It may be shown that its extension in the complex plan, $z \in \mathbb{C}$, is a meromorphic function. The poles of the gamma function are any negative integer values of $z$. Then for any $k \in \mathbb{Z}_-$

$$\lim_{z \to k} \Gamma(z) = \infty \qquad (2.4)$$

Some properties of the gamma function may be demonstrated for any $z \in \mathbb{C} - \{k/k \in \mathbb{Z}_-\}$. We have the multiplication theorem

$$\Gamma(z)\Gamma\left(z + \frac{1}{m}\right)\Gamma\left(z + \frac{2}{m}\right)\ldots\Gamma\left(z + \frac{m-1}{m}\right) =$$
$$= (2\pi)^{\frac{m-1}{2}} m^{\frac{1}{2} - mz} \Gamma(mz) \qquad (2.5)$$

particularly

$$\Gamma(z)\Gamma(z+\frac{1}{2}) = (2)^{1-2z}\sqrt{\pi}\Gamma(2z) \qquad (2.6)$$

For $n \in \mathbb{N}$, we have

$$\Gamma\left(\frac{1}{2}+n\right) = \frac{(2n)!}{4^n n!}\sqrt{\pi} \quad \Gamma\left(\frac{1}{2}-n\right) = \frac{(-4)^n n!}{(2n)!}\sqrt{\pi} \quad (2.7)$$

We give some useful numerical values of the gamma functions

$$\begin{cases} \Gamma\left(-\frac{5}{2}\right) = -\frac{8\sqrt{\pi}}{15} \\ \Gamma\left(-\frac{3}{2}\right) = \frac{4\sqrt{\pi}}{3} \quad \Gamma\left(-\frac{1}{2}\right) = -2\sqrt{\pi} \\ \Gamma\left(\frac{1}{2}\right) = \sqrt{\pi} \quad \Gamma\left(\frac{3}{2}\right) = \frac{\sqrt{\pi}}{2} \\ \Gamma\left(\frac{5}{2}\right) = \frac{3\sqrt{\pi}}{4} \end{cases} \qquad (2.8)$$

# 3. First Definition

## 3.1. First Definition

For any doublet $(s, \alpha) \in \mathbb{Q} \times \mathbb{R} - \{(-1,-1)\}$ and for $x \in \mathbb{R}^+$, a derivative of order $s$ of the function $f$ such

$$f(x) = ax^\alpha = \begin{cases} a & \text{if } x = 0 \text{ and } \alpha = 0 \\ 0 & \text{if } x = 0 \text{ and } \alpha \neq 0 \\ ae^{\alpha \ln(x)} & \text{if } x > 0 \text{ and } \alpha \neq 0 \end{cases} \quad (3.1.1)$$

is defined by

$$\frac{d^s(ax^\alpha)}{dx^s} = a\frac{d^s x^\alpha}{dx^s}$$

$$= a\frac{\Gamma(\alpha+1)}{\Gamma(\alpha-s+1)}x^{\alpha-s} + \sum_{k=-\infty}^{-1} c_k \frac{x^{k-s}}{\Gamma(k-s+1)} \qquad (3.1.2)$$

in which $\{c_k\}$ is a set of arbitrary constants with finite values. We will see that the second part of the relation (3.1.2) is the terms which allow the unified expression of fractional derivative to be applicable for both negative and positive orders. In the case of positive integer orders, this second term of (3.1.2) is equal to zero. In the case of negative integer orders, this term becomes a finite summation. So there is no problem about the convergence in the case of integer order. But in the case of fractional orders, this term is an infinite summation and we have to discuss about the convergence. We will see also that because of the existence of the arbitrary constants $c_k$, the fractional derivative is not unique when the order is not a positive integer.

***Theorem 3.1***

For $s = n \in \mathbb{N}$, the derivative of order $s$ as defined in (3.1.2) is equivalent to the ordinary derivative

Proof

For $s = n \in \mathbb{N}$, we have from the relation (3.1.2)

$$\frac{d^n(ax^\alpha)}{dx^n} = a\frac{d^s x^\alpha}{dx^s}$$

$$= a\frac{\Gamma(\alpha+1)}{\Gamma(\alpha-n+1)}x^{\alpha-n} + \sum_{k=-\infty}^{-1} c_k \frac{x^{k-n}}{\Gamma(k-n+1)} \qquad (3.1.3)$$

Because of the property (2.4), the second part of the expression (3.1.3) is equal to zero for positive integer $n$. Then

$$\frac{d^n(ax^\alpha)}{dx^n} = a\frac{\Gamma(\alpha+1)}{\Gamma(\alpha-n+1)}x^{\alpha-n}$$

$$= a\alpha(\alpha-1)\ldots(\alpha-n)x^{\alpha-n} \qquad (3.1.4)$$

The relation (3.14) is the expression of ordinary derivative of order $n$.

Example

$$\frac{d^6 x^7}{dx^6} = \frac{\Gamma(7+1)}{\Gamma(7-6+1)}x^{7-6} = \frac{7!}{1!}x = 5040x$$

$$\frac{d^3(3x^6)}{dx^3} = 3\frac{\Gamma(6+1)}{\Gamma(6-3+1)}x^{6-3} = 3\frac{6!}{3!}x^3 = 360x^3$$

$$\frac{d^5(-8x^{11})}{dx^5} = -8\frac{11!}{6!}x^6 = -443520x^6$$

***Theorem 3.2***

If the order is a negative integer $s = -n$ with $n \in \mathbb{N}$, the fractional derivative as defined in (3.1.2) is not unique. It depends on $n$ arbitrary constants.

$$\frac{d^{-n}(ax^\alpha)}{dx^{-n}} = a\frac{d^{-n} x^\alpha}{dx^{-n}}$$

$$= a\frac{\Gamma(\alpha+1)}{\Gamma(\alpha-n+1)}x^{\alpha+n} + \sum_{k=-n}^{-1} c_k \frac{x^{k+n}}{(k+n)!} \qquad (3.1.5)$$

Proof

If $s = -n$ with $n \in \mathbb{N}$ and according to the properties (2.3) and (2.4) of the gamma function. We have

$$\Gamma(k-s+1) = \Gamma(k+n+1)$$

$$= \begin{cases} \infty & \text{if } k < -n \\ (k+n)! & \text{if } k \geq -n \end{cases} \qquad (3.1.6)$$

From the relation (3.1.2) and the relation (3.1.6) we may deduce easily the relation (3.1.5).

$$\frac{d^{-n}(ax^\alpha)}{dx^{-n}} = a\frac{\Gamma(\alpha+1)}{\Gamma(\alpha-n+1)}x^{\alpha+n}$$

$$+ \sum_{k=-n}^{-1} c_k \frac{x^{k+n}}{(k+n)!} \qquad (3.1.7)$$

The right hand side is the $n$-times indefinite integral (primitive) of $ax^\alpha$. If we derivate $n$ times the relation

(3.1.7) all the terms containing the constants $c_k$ disappear.
Examples

$$\frac{d^{-1}x^3}{dx^{-1}} = \frac{\Gamma(3+1)}{\Gamma(3+1+1)}x^{3+1} + \sum_{k=-1}^{-1} c_k \frac{x^{k+1}}{\Gamma(k+1+1)}$$

$$= \frac{x^4}{4} + c_{-1}$$

$$\frac{d^{-3}5x}{dx^{-3}} = 5\frac{\Gamma(1+1)}{\Gamma(1+3+1)}x^{1+3} + \sum_{k=-3}^{-1} c_k \frac{x^{k+3}}{\Gamma(k+3+1)}$$

$$= 5\frac{x^4}{4} + c_{-1}\frac{x^2}{2} + c_{-2}x + c_{-3}$$

**Theorem 3.3**

If the order $s$ is a negative integer, $s = -n$ with $n \in \mathbb{N}$, the fractional derivative as defined by the relation (3.1.2) is not unique. It depends on $n$ arbitrary constants.

Proof

It is evident by looking at the second part of the relation (3.1.7). If $F(x)$ is a derivative of order $s$ of $f(x) = ax^\alpha$, then the function

$$G(x) = F(x) + \sum_{k=-n}^{-1} b_k \frac{x^{k+n}}{(k+n)!} \quad (3.1.8)$$

is also a derivative of order $s$ of $f(x)$ in which $\{b_k\}$ is a set of arbitrary finite constant numbers. We may introduce the notion of semi-equality because of the relation (3.1.8) between $G(x)$ and $F(x)$. This relation is an equivalence relation.

**Theorem 3.4**

If the order $s$ is a fraction $s = \frac{p}{q} \in \mathbb{Q}$, a sufficient condition for the fractional derivative, as defined in (3.1.2) to be well defined for any $x \in \mathbb{R}_+$ is

$$\lim_{k \to -\infty} \left| k \frac{c_{k-1}}{c_k} \right| = 0 \quad (3.1.9)$$

Proof
The relation (3.1.2) gives for $s = \frac{p}{q}$

$$\frac{d^s(ax^\alpha)}{dx^s} = \frac{d^{\frac{p}{q}}(ax^\alpha)}{dx^{\frac{p}{q}}}$$

$$= a\frac{\Gamma(\alpha+1)}{\Gamma(\alpha-\frac{p}{q}+1)}x^{\alpha-\frac{p}{q}} + \sum_{k=-\infty}^{-1} c_k \frac{x^{k-\frac{p}{q}}}{\Gamma(k-\frac{p}{q}+1)} \quad (3.1.10)$$

The last members of this relation is well defined if the infinite summation is convergent. Let be

$$R(s,x) = \sum_{k=-\infty}^{-1} c_k \frac{x^{k-s}}{\Gamma(k-s+1)} = \sum_{n=1}^{+\infty} c_{-n} \frac{x^{-n-s}}{\Gamma(-n-s+1)}$$

and

$$r_n(s,x) = c_{-n} \frac{x^{-n-s}}{\Gamma(-n-s+1)}$$

Then

$$R(s,x) = \sum_{n=1}^{+\infty} r_n(s,x)$$

According to the Abel's criterium, the series $R(s,x)$ is convergent if and only if

$$\lim_{n \to +\infty} \left| \frac{r_{n+1}(s,x)}{r_n(s,x)} \right| < 1 \Leftrightarrow \lim_{n \to +\infty} \left| -(n+s)\frac{c_{-(n+1)}}{c_{-n}x} \right| < 1$$

$$\Leftrightarrow |x| > \lim_{n \to +\infty} \left| -(n+s)\frac{c_{-(n+1)}}{c_{-n}} \right|$$

So the definition (3.1.10) is meaningful for all $x > 0$ if the coefficients $c_{-n} = c_k$ fulfill the condition

$$\lim_{k \to -\infty} \left| (k-s)\frac{c_{k-1}}{c_k} \right| = 0 \Leftrightarrow \lim_{k \to -\infty} \left| k\frac{c_{k-1}}{c_k} \right| = 0$$

In all following calculations, expressions and use of the above defined fractional derivative, we assume that the condition (3.1.9) is fulfilled. This condition gives a restriction on the possible values of constants $c_k$.

**Theorem 3.5**

If the order $s$ is a fraction $s = \frac{p}{q} \in \mathbb{Q}$, the fractional derivative, as defined in (3.1.2) is not unique.

**Proof**

The relation (3.1.2) gives for $s = \frac{p}{q}$ the relation (3.1.10). According to this relation, we remark that if $F(x)$ is a derivative of order $\frac{p}{q}$ of the function $f(x) = ax^\alpha$, then the function

$$G(x) = F(x) + \sum_{k=-\infty}^{-1} b_k \frac{x^{k-\frac{p}{q}}}{\Gamma(k-\frac{p}{q}+1)} \quad (3.1.11)$$

is also a derivative of order $\frac{p}{q}$ of $f(x)$. The constants $b_k$ are arbitrary constants but they have to fulfill the conditions (3.1.9) for the convergence of the infinite summation. For the relation (3.1.11), we recall and extend the notion of "semi-equality" as defined in the relation (3.1.8).

### 3.2. Properties

#### 3.2.1 Linearity and Semi-Linearity

**Theorem 3.6**

If the order $s$ is positive integer, the fractional derivative of order $s$ verifies the property of linearity

$$\frac{d^s(ax^\alpha + bx^\beta)}{dx^s} = a\frac{d^s x^\alpha}{dx^s} + b\frac{d^s x^\beta}{dx^s} \quad (3.2.1)$$

Proof

This property is evident according to the theorem 3.1. In the case of positive integer order, the fractional derivative is equivalent to the ordinary derivative which is well known to be linear at any order.

**Theorem 3.7**

If the order $s$ is a negative integer, or a negative or positive fraction, the fractional derivative of order $s$ does not have the property of linearity. But it has a property that we may define as "semi-linearity"

$$\frac{d^s(ax^\alpha + bx^\beta)}{dx^s}$$

$$= a\frac{d^s x^\alpha}{dx^s} + b\frac{d^s x^\beta}{dx^s} + \sum_{k=-\infty}^{-1} c_k \frac{x^{k-s}}{\Gamma(k-s+1)} \quad (3.2.2)$$

Proof

If $s$ is not a positive integer, we have

$$a\frac{d^s x^\alpha}{dx^s} + b\frac{d^s x^\beta}{dx^s}$$

$$= a\frac{\Gamma(\alpha+1)}{\Gamma(\alpha-s+1)}x^{\alpha-s} + \sum_{k=-\infty}^{-1} c_{k1}\frac{x^{k-s}}{\Gamma(k-s+1)}$$

$$+ b\frac{\Gamma(\beta+1)}{\Gamma(\beta-s+1)}x^{\beta-s} + \sum_{k=-\infty}^{-1} c_{k2}\frac{x^{k-s}}{\Gamma(k-s+1)} \quad (3.2.3)$$

$$\frac{d^s(ax^\alpha + bx^\beta)}{dx^s}$$

$$= a\frac{\Gamma(\alpha+1)}{\Gamma(\alpha-s+1)}x^{\alpha-s} + b\frac{\Gamma(\beta+1)}{\Gamma(\beta-s+1)}x^{\beta-s}$$

$$+ \sum_{k=-\infty}^{-1} c_{k3}\frac{x^{k-s}}{\Gamma(k-s+1)} \quad (3.2.4)$$

Because of the arbitrariness on the values of the constants $c_{k1}, c_{k2}$, and $c_k$, we may identify $c_k = c_{k3} - c_{k1} - c_{k2}$ and then deduce from the relation (3.2.3) and (3.2.4) the relation (3.2.2). We may qualify the property (3.2.2) as "semi-linearity" on the same footing as semi-equality defined in the relations (3.1.8) and (3.1.11).

### 3.2.2. Commutativity and Semi-Commutativity

**Theorem 3.8**

If the order is a positive integer, the fractional derivative has the property of commutativity

$$\frac{d^{s_1}}{dx^{s_1}}\frac{d^{s_2}(ax^\alpha)}{dx^{s_2}} = \frac{d^{s_2}}{dx^{s_2}}\frac{d^{s_1}(ax^\alpha)}{dx^{s_1}}$$

$$= \frac{d^{s_1+s_2}(ax^\alpha)}{dx^{s_1+s_2}} \quad (3.2.5)$$

Proof

This property is evident because of the Theorem 3.1, in the case of positive integer order, the fractional derivative is exactly the ordinary derivative which is well known to have the property of commutativity.

**Theorem 3.9**

If the order $s$ is a negative integer, or a negative or positive fraction, the fractional derivative of order $s$ does not have the property of commutativity. But it has a property that we may define as "semi-commutativity" (semi-group property)

Proof

using the definition (3.1.2) we obtain

$$\frac{d^{s_2}}{dx^{s_2}}\frac{d^{s_1}(ax^\alpha)}{dx^{s_1}} = a\frac{\Gamma(\alpha+1)}{\Gamma(\alpha-(s_1+s_2)+1)}x^{\alpha-(s_1+s_2)}$$

$$+ \sum_{k=-\infty}^{-1} c_k^{21}\frac{x^{k-(s_1+s_2)}}{\Gamma(k-s_1+s_2+1)} + \sum_{k=-\infty}^{-1} c_k^{2}\frac{x^{k-s_2}}{\Gamma(k-s+1)}$$

$$= \frac{d^{s_1+s_2}(ax^\alpha)}{dx^{s_1+s_2}} + \sum_{k=-\infty}^{-1} c_k^{2}\frac{x^{k-s_2}}{\Gamma(k-s+1)} \quad (3.2.6)$$

$$\frac{d^{s_1}}{dx^{s_1}}\frac{d^{s_2}(ax^\alpha)}{dx^{s_2}} = a\frac{\Gamma(\alpha+1)}{\Gamma(\alpha-(s_1+s_2)+1)}x^{\alpha-(s_1+s_2)}$$

$$+ \sum_{k=-\infty}^{-1} c_k^{12}\frac{x^{k-(s_1+s_2)}}{\Gamma(k-s_1+s_2+1)} + \sum_{k=-\infty}^{-1} c_k^{1}\frac{x^{k-s_1}}{\Gamma(k-s+1)}$$

$$= \frac{d^{s_1+s_2}(ax^\alpha)}{dx^{s_1+s_2}} + \sum_{k=-\infty}^{-1} c_k^{1}\frac{x^{k-s_1}}{\Gamma(k-s+1)} \quad (3.2.7)$$

Then we have

$$\frac{d^{s_1}}{dx^{s_1}}\frac{d^{s_2}(ax^\alpha)}{dx^{s_2}} - \frac{d^{s_2}}{dx^{s_2}}\frac{d^{s_1}(ax^\alpha)}{dx^{s_1}} =$$

$$\sum_{k=-\infty}^{-1} c_k^{2}\frac{x^{k-s_2}}{\Gamma(k-s+1)} - \sum_{k=-\infty}^{-1} c_k^{1}\frac{x^{k-s_1}}{\Gamma(k-s+1)} \quad (3.2.8)$$

According to the relation (3.2.8) we don't have in general the property of commutativity. But on equal footing as the definitions of "semi-equality" and 'semi-linearity" defined in the relations (3.1.8), (3.1.11) and (3.2.2),we may define the property of "semi-commutativity" for the relation (3.2.8).

### 3.2.3. Important Relation Verified by the Fractional De-

*rivatives.*

For any $\frac{p}{q} \in \mathbb{Q}$, we may always assume, without loss of generality, that $p \in \mathbb{Z}$ and $q \in \mathbb{N}^*$. We have the theorem:

**Theorem 3.9**

For $n \in \mathbb{N}$ and $n \leq q$

$$\underbrace{\frac{d^{\frac{p}{q}}}{dx^{\frac{p}{q}}} \cdot \frac{d^{\frac{p}{q}}}{dx^{\frac{p}{q}}} \cdots \frac{d^{\frac{p}{q}}}{dx^{\frac{p}{q}}}}_{n \text{ times}} (ax^\alpha)$$

$$= \frac{d^{\frac{np}{q}}(ax^\alpha)}{dx^{\frac{np}{q}}} + \sum_{m=0}^{n-1} \sum_{k=-\infty}^{-1} c_k^m \frac{x^{k-\frac{mp}{q}}}{\Gamma(k-\frac{mp}{q}+1)} \quad (3.2.9)$$

Proof

We use a mathematical induction.

For $q = 1$, the equality is evident.

For $q = 2$ and $n = 1$ according to (3.1.2) we have

$$\frac{d^{\frac{p}{2}}(ax^\alpha)}{dx^{\frac{p}{2}}} = a \frac{\Gamma(\alpha+1)}{\Gamma\left(\alpha-\frac{p}{2}+1\right)} x^{\alpha-\frac{p}{2}}$$

$$+ \sum_{k=-\infty}^{-1} c_k \frac{x^{k-\frac{p}{2}}}{\Gamma(k-\frac{p}{2}+1)}$$

For $q = 2$ and $n = 2$

$$\frac{d^{\frac{p}{2}}}{dx^{\frac{p}{2}}} \left[\frac{d^{\frac{p}{2}}(ax^\alpha)}{dx^{\frac{p}{2}}}\right] = \frac{d^{\frac{p}{2}}}{dx^{\frac{p}{2}}} \left[a \frac{\Gamma(\alpha+1)}{\Gamma\left(\alpha-\frac{p}{2}+1\right)} x^{\alpha-\frac{p}{2}}\right.$$

$$\left. + \sum_{k=-\infty}^{-1} c_k \frac{x^{k-\frac{p}{2}}}{\Gamma(k-\frac{p}{2}+1)}\right]$$

$$= a \frac{\Gamma(\alpha+1)}{\Gamma\left(\alpha-\frac{p}{2}+1\right)} \cdot \frac{\Gamma\left(\alpha-\frac{p}{2}+1\right)}{\Gamma\left(\alpha-\frac{p}{2}-\frac{p}{2}+1\right)} x^{\alpha-\frac{p}{2}-\frac{p}{2}}$$

$$+ \sum_{k=-\infty}^{-1} c_k \frac{1}{\Gamma(k-\frac{p}{2}+1)} \frac{\Gamma\left(k-\frac{p}{2}+1\right)}{\Gamma\left(k-\frac{p}{2}-\frac{p}{2}+1\right)} x^{k-\frac{p}{2}-\frac{p}{2}}$$

$$+ \sum_{k=-\infty}^{-1} c_k^1 \frac{x^{k-\frac{p}{2}}}{\Gamma(k-\frac{p}{2}+1)}$$

$$= \frac{d^p(ax^\alpha)}{dx^p} + \sum_{k=-\infty}^{-1} c_k^1 \frac{x^{k-\frac{p}{2}}}{\Gamma(k-\frac{p}{2}+1)} \quad (3.2.10)$$

For $q = 3$ and $n = 1$

$$\frac{d^{\frac{p}{3}}(ax^\alpha)}{dx^{\frac{p}{3}}} = a \frac{\Gamma(\alpha+1)}{\Gamma\left(\alpha-\frac{p}{3}+1\right)} x^{\alpha-\frac{p}{3}}$$

$$+ \sum_{k=-\infty}^{-1} c_k \frac{x^{k-\frac{p}{3}}}{\Gamma(k-\frac{p}{3}+1)} \quad (3.2.11)$$

For $q = 3$ and $n = 2$

$$\frac{d^{\frac{p}{3}}}{dx^{\frac{p}{3}}} \left[\frac{d^{\frac{p}{3}}(ax^\alpha)}{dx^{\frac{p}{3}}}\right] = \frac{d^{\frac{p}{3}}}{dx^{\frac{p}{3}}} \left[a \frac{\Gamma(\alpha+1)}{\Gamma\left(\alpha-\frac{p}{3}+1\right)} x^{\alpha-\frac{p}{3}}\right.$$

$$\left. + \sum_{k=-\infty}^{-1} c_k \frac{x^{k-\frac{p}{3}}}{\Gamma(k-\frac{p}{3}+1)}\right]$$

$$= a \frac{\Gamma(\alpha+1)}{\Gamma\left(\alpha-\frac{p}{3}+1\right)} \cdot \frac{\Gamma\left(\alpha-\frac{p}{3}+1\right)}{\Gamma\left(\alpha-\frac{p}{3}-\frac{p}{3}+1\right)} x^{\alpha-\frac{p}{3}-\frac{p}{3}}$$

$$+ \sum_{k=-\infty}^{-1} c_k \frac{1}{\Gamma(k-\frac{p}{3}+1)} \frac{\Gamma\left(k-\frac{p}{3}+1\right)}{\Gamma\left(k-\frac{p}{3}-\frac{p}{3}+1\right)} x^{k-\frac{p}{2}-\frac{p}{2}}$$

$$+ \sum_{k=-\infty}^{-1} c_k^1 \frac{x^{k-\frac{p}{3}}}{\Gamma(k-\frac{p}{3}+1)}$$

$$= \frac{d^{\frac{2p}{3}}(ax^\alpha)}{dx^{\frac{2p}{3}}} + \sum_{k=-\infty}^{-1} c_k^1 \frac{x^{k-\frac{p}{3}}}{\Gamma(k-\frac{p}{3}+1)}$$

$$= \frac{d^{\frac{2p}{3}}(ax^\alpha)}{dx^{\frac{2p}{3}}} + \sum_{m=0}^{1} \sum_{k=-\infty}^{-1} c_k^m \frac{x^{k-\frac{mp}{3}}}{\Gamma(k-\frac{mp}{3}+1)} \quad (3.2.12)$$

For $q = 3$ and $n = 3$

$$\frac{d^{\frac{p}{3}}}{dx^{\frac{p}{3}}} \frac{d^{\frac{p}{3}}}{dx^{\frac{p}{3}}} \frac{d^{\frac{p}{3}}(ax^\alpha)}{dx^{\frac{p}{3}}} = \frac{d^{\frac{p}{3}}}{dx^{\frac{p}{3}}} \left[a \frac{\Gamma(\alpha+1)}{\Gamma\left(\alpha-\frac{2p}{3}+1\right)} x^{\alpha-\frac{2p}{3}}\right.$$

$$+ \sum_{k=-\infty}^{-1} c_k \frac{x^{k-\frac{2p}{3}}}{\Gamma(k-\frac{2p}{3}+1)} + \sum_{k=-\infty}^{-1} c_k^1 \frac{x^{k-\frac{p}{3}}}{\Gamma(k-\frac{p}{3}+1)}\right]$$

$$= a \frac{\Gamma(\alpha+1)}{\Gamma\left(\alpha-\frac{2p}{3}+1\right)} \cdot \frac{\Gamma\left(\alpha-\frac{2p}{3}+1\right)}{\Gamma\left(\alpha-\frac{2p}{3}-\frac{p}{3}+1\right)} x^{\alpha-\frac{2p}{3}-\frac{p}{3}}$$

$$+ \sum_{k=-\infty}^{-1} c_k \frac{\Gamma\left(k-\frac{2p}{3}+1\right)}{\Gamma\left(k-\frac{2p}{3}-\frac{p}{3}+1\right)\Gamma(k-\frac{2p}{3}+1)} x^{k-\frac{2p}{3}-\frac{p}{3}}$$

$$+ \sum_{k=-\infty}^{-1} c_k^1 \frac{\Gamma\left(\alpha-\frac{p}{3}+1\right)}{\Gamma\left(\alpha-\frac{p}{3}-\frac{p}{3}+1\right)\Gamma\left(k-\frac{p}{3}+1\right)} x^{k-\frac{p}{3}-\frac{p}{3}}$$

$$+ \sum_{k=-\infty}^{-1} c_k^2 \frac{x^{k-\frac{p}{3}}}{\Gamma(k-\frac{p}{3}+1)}$$

$$= \frac{d^p(ax^\alpha)}{dx^p} + \sum_{m=0}^{2}\sum_{k=-\infty}^{-1} c_k^m \frac{x^{k-\frac{mp}{3}}}{\Gamma(k-\frac{mp}{3}+1)} \quad (3.2.13)$$

For any $n \leq q, q \in \mathbb{N}^*$, we assume the hypothesis (3.2.9)

$$\underbrace{\frac{d^{\frac{p}{q}}}{dx^{\frac{p}{q}}} \cdot \frac{d^{\frac{p}{q}}}{dx^{\frac{p}{q}}} \cdots \frac{d^{\frac{p}{q}}}{dx^{\frac{p}{q}}}}_{n \text{ times}}(ax^\alpha)$$

$$= \frac{d^{\frac{np}{q}}(ax^\alpha)}{dx^{\frac{np}{q}}} + \sum_{m=0}^{n-1}\sum_{k=-\infty}^{-1} c_k^m \frac{x^{k-\frac{mp}{q}}}{\Gamma(k-\frac{mp}{q}+1)} \quad (3.2.14)$$

We have now to prove that this formula remains true for $n+1 < q$ i.e

$$\underbrace{\frac{d^{\frac{p}{q}}}{dx^{\frac{p}{q}}} \cdot \frac{d^{\frac{p}{q}}}{dx^{\frac{p}{q}}} \cdots \frac{d^{\frac{p}{q}}}{dx^{\frac{p}{q}}}}_{(n+1) \text{ times}}(ax^\alpha)$$

$$= \frac{d^{\frac{(n+1)p}{q}}(ax^\alpha)}{dx^{\frac{(n+1)p}{q}}} + \sum_{m=0}^{n}\sum_{k=-\infty}^{-1} c_k^m \frac{x^{k-\frac{mp}{q}}}{\Gamma(k-\frac{mp}{q}+1)} \quad (3.2.15)$$

Let us perform the calculation

$$\underbrace{\frac{d^{\frac{p}{q}}}{dx^{\frac{p}{q}}} \cdot \frac{d^{\frac{p}{q}}}{dx^{\frac{p}{q}}} \cdots \frac{d^{\frac{p}{q}}}{dx^{\frac{p}{q}}}}_{(n+1) \text{ times}}(ax^\alpha) = \frac{d^{\frac{p}{q}}}{dx^{\frac{p}{q}}}[\underbrace{\frac{d^{\frac{p}{q}}}{dx^{\frac{p}{q}}} \cdot \frac{d^{\frac{p}{q}}}{dx^{\frac{p}{q}}} \cdots \frac{d^{\frac{p}{q}}}{dx^{\frac{p}{q}}}}_{n \text{ times}}(ax^\alpha)]$$

$$= \frac{d^{\frac{p}{q}}}{dx^{\frac{p}{q}}}\left[\frac{d^{\frac{np}{q}}(ax^\alpha)}{dx^{\frac{np}{q}}} + \sum_{m=0}^{n-1}\sum_{k=-\infty}^{-1} c_k^m \frac{x^{k-\frac{mp}{q}}}{\Gamma\left(k-\frac{mp}{q}+1\right)}\right]$$

$$= \frac{d^{\frac{p}{q}}}{dx^{\frac{p}{q}}}\left[a\frac{\Gamma(\alpha+1)}{\Gamma(\alpha-\frac{np}{q}+1)}x^{\alpha-\frac{np}{q}} + \sum_{k=-\infty}^{-1} c_k \frac{x^{k-\frac{np}{q}}}{\Gamma(k-\frac{np}{q}+1)}\right.$$

$$\left.+ \sum_{m=0}^{n-1}\sum_{k=-\infty}^{-1} c_k^m \frac{x^{k-\frac{mp}{q}}}{\Gamma\left(k-\frac{mp}{q}+1\right)}\right]$$

$$= a\frac{\Gamma(\alpha+1)}{\Gamma(\alpha-\frac{np}{q}+1)}\frac{\Gamma(\alpha-\frac{np}{q}+1)}{\Gamma(\alpha-\frac{np}{q}-\frac{p}{q}+1)}x^{\alpha-\frac{np}{q}-\frac{p}{q}}$$

$$+ \sum_{k=-\infty}^{-1} c_k \frac{\Gamma(k-\frac{np}{q}+1)}{\Gamma(k-\frac{np}{q}-\frac{p}{q}+1)}\frac{x^{k-\frac{np}{q}}}{\Gamma(k-\frac{np}{q}+1)}$$

$$+ \sum_{m=0}^{n-1}\sum_{k=-\infty}^{-1} c_k^m \frac{\Gamma(k-\frac{mp}{q}+1)}{\Gamma(k-\frac{mp}{q}-\frac{p}{q}+1)}\frac{x^{k-\frac{mp}{q}}}{\Gamma\left(k-\frac{mp}{q}+1\right)}$$

$$+ \sum_{k=-\infty}^{-1} c_k^n \frac{x^{k-\frac{p}{q}}}{\Gamma(k-\frac{np}{q}+1)}$$

$$= a\frac{\Gamma(\alpha+1)}{\Gamma(\alpha-\frac{(n+1)p}{q}+1)}x^{\alpha-\frac{(n+1)p}{q}}$$

$$+ \sum_{k=-\infty}^{-1} c_k \frac{x^{k-\frac{np}{q}}}{\Gamma\left(k-\frac{(n+1)p}{q}+1\right)}$$

$$+ \sum_{m=1}^{n-1}\sum_{k=-\infty}^{-1} c_k^m \frac{x^{k-\frac{mp}{q}}}{\Gamma\left(k-\frac{mp}{q}-\frac{p}{q}+1\right)}$$

$$+ \sum_{k=-\infty}^{-1} c_k^n \frac{x^{k-\frac{p}{q}}}{\Gamma(k-\frac{np}{q}+1)}$$

$$= \frac{d^{\frac{(n+1)p}{q}}(ax^\alpha)}{dx^{\frac{(n+1)p}{q}}} + \sum_{m=0}^{n}\sum_{k=-\infty}^{-1} c_k^m \frac{x^{k-\frac{mp}{q}}}{\Gamma(k-\frac{mp}{q}+1)}$$

## 4. Second Definition

Replacing the order $s$ by its entire part in the summation which appears in the relation (3.1.2), we may give another definition of the fractional derivative which looks simpler because it eliminates the complexity introduced by the existence of the infinite summation. So we do not have to study the convergence of this term.

### 4.1. Second Definition

For any doublet $(s, \alpha) \in \mathbb{Q} \times \mathbb{R} - \{(-1, -1)\}$ and for $x \in \mathbb{R}^+$, a derivative of order $s$ of the function $f$ such

$$f(x) = ax^\alpha = \begin{cases} a & \text{if } x = 0 \text{ and } \alpha = 0 \\ 0 & \text{if } x = 0 \text{ and } \alpha \neq 0 \\ ae^{\alpha \ln(x)} & \text{if } x > 0 \text{ and } \alpha \neq 0 \end{cases} \quad (4.1.1)$$

is defined by

$$\frac{d^s(ax^\alpha)}{dx^s} = a\frac{d^s x^\alpha}{dx^s}$$

$$= a\frac{\Gamma(\alpha+1)}{\Gamma(\alpha-s+1)}x^{\alpha-s} + \sum_{k=-|[s]|}^{-1} c_k \frac{x^{k-s}}{\Gamma(k-[s]+1)} \quad (4.1.2)$$

$[s]$ is the entire part of $s$. This second definition gives exactly the same results as the definition (3.1.2) in the case of integer order. But some differences appear between the two definitions in the case of fractional order.

### Theorem 4.1

For $s = n$ with $n \in \mathbb{N}$, the derivative of order $s$ as defined in (4.1.2) is the ordinary derivative of order $n$.

Proof
From the definition (4.1.2) and using the property (2.4) of the gamma function, we obtain

$$\frac{d^n(ax^\alpha)}{dx^n} = a\frac{\Gamma(\alpha+1)}{\Gamma(\alpha-n+1)}x^{\alpha-n}$$

$$= a\alpha(\alpha-1)\ldots(\alpha-n)x^{\alpha-n} \quad (4.1.3)$$

The relation (4.1.3) is the ordinary expression of derivative of order $n$ for $n \in \mathbb{N}$.

**Theorem 4.2**

For $s = \frac{p}{q} \in \mathbb{Q}^+$, the derivative of order $s$ as defined in (4.1.2) is given by the expression

$$\frac{d^{\frac{p}{q}}(ax^\alpha)}{dx^{\frac{p}{q}}} = a\frac{d^{\frac{p}{q}}(x^\alpha)}{dx^{\frac{p}{q}}} = a\frac{\Gamma(\alpha+1)}{\Gamma(\alpha-\frac{p}{q}+1)}x^{\alpha-\frac{p}{q}} \quad (4.1.4)$$

Proof
The relation (4.1.4) is deduced easily from the relation (4.1.2) by replacing $s = \frac{p}{q}$ and noting that the second part of this expression is equal to 0 because of the relation (2.4). According to (4.1.3) and (4.1.4), we remark that for positive value of the order $s$ (integer or fractional), the fractional derivative as defined in (4.1.2) is unique. We will see that it is not always the case when $s$ is negative.

**Theorem 4.3**

For $s = \frac{p}{q} \in \mathbb{Q}^+$, with $p, n \in \mathbb{N}, q \in \mathbb{N}^*$ we have the equality

$$\underbrace{\frac{d^{\frac{p}{q}}}{dx^{\frac{p}{q}}} \cdot \frac{d^{\frac{p}{q}}}{dx^{\frac{p}{q}}} \ldots \frac{d^{\frac{p}{q}}}{dx^{\frac{p}{q}}}}_{n\ times}(ax^\alpha) = \frac{d^{\frac{np}{q}}(ax^\alpha)}{dx^{\frac{np}{q}}} \quad (4.1.5)$$

and particularly

$$\underbrace{\frac{d^{\frac{p}{q}}}{dx^{\frac{p}{q}}} \cdot \frac{d^{\frac{p}{q}}}{dx^{\frac{p}{q}}} \ldots \frac{d^{\frac{p}{q}}}{dx^{\frac{p}{q}}}}_{q\ times}(ax^\alpha) = \frac{d^p(ax^\alpha)}{dx^p} \quad (4.1.6)$$

Proof
We utilize mathematical induction.
For $q = 1$ the relation is evident
For $q = 2$ and $n = 1$, according to (4.1.4), we have

$$\frac{d^{\frac{1}{2}}(ax^\alpha)}{dx^{\frac{1}{2}}} = a\frac{d^{\frac{1}{2}}(x^\alpha)}{dx^{\frac{1}{2}}} = a\frac{\Gamma(\alpha+1)}{\Gamma(\alpha-\frac{p}{q}+1)}x^{\alpha-\frac{1}{2}} \quad (4.1.7)$$

For $q = 2$ and $n = 2$, according to (4.14), we have

$$\frac{d^{\frac{1}{2}}}{dx^{\frac{1}{2}}}\frac{d^{\frac{1}{2}}(x^\alpha)}{dx^{\frac{1}{2}}} = \frac{d^{\frac{1}{2}}}{dx^{\frac{1}{2}}}\left[a\frac{\Gamma(\alpha+1)}{\Gamma(\alpha-\frac{1}{2}+1)}x^{\alpha-\frac{1}{2}}\right]$$

$$= a\frac{\Gamma(\alpha+1)}{\Gamma(\alpha-1+1)} = a\alpha x^{\alpha-1} = \frac{dax^\alpha}{dx} \quad (4.1.8)$$

For any $n \in \mathbb{N}, q \in \mathbb{N}^*$, we assume that

$$\underbrace{\frac{d^{\frac{p}{q}}}{dx^{\frac{p}{q}}} \cdot \frac{d^{\frac{p}{q}}}{dx^{\frac{p}{q}}} \ldots \frac{d^{\frac{p}{q}}}{dx^{\frac{p}{q}}}}_{n\ times}(ax^\alpha) = \frac{d^{\frac{np}{q}}(ax^\alpha)}{dx^{\frac{np}{q}}} \quad (4.1.9)$$

We have now to prove that this formula remains true for $n + 1$

$$\underbrace{\frac{d^{\frac{p}{q}}}{dx^{\frac{p}{q}}} \cdot \frac{d^{\frac{p}{q}}}{dx^{\frac{p}{q}}} \ldots \frac{d^{\frac{p}{q}}}{dx^{\frac{p}{q}}}}_{(n+1)\ times}(ax^\alpha) = \frac{d^{\frac{(n+1)p}{q}}(ax^\alpha)}{dx^{\frac{(n+1)p}{q}}} \quad (4.1.10)$$

Let us perform the calculation

$$\underbrace{\frac{d^{\frac{p}{q}}}{dx^{\frac{p}{q}}} \cdot \frac{d^{\frac{p}{q}}}{dx^{\frac{p}{q}}} \ldots \frac{d^{\frac{p}{q}}}{dx^{\frac{p}{q}}}}_{(n+1)\ times}(ax^\alpha) = \frac{d^{\frac{p}{q}}}{dx^{\frac{p}{q}}}[\underbrace{\frac{d^{\frac{p}{q}}}{dx^{\frac{p}{q}}} \cdot \frac{d^{\frac{p}{q}}}{dx^{\frac{p}{q}}} \ldots \frac{d^{\frac{p}{q}}}{dx^{\frac{p}{q}}}}_{n\ times}(ax^\alpha)]$$

$$= \frac{d^{\frac{p}{q}}}{dx^{\frac{p}{q}}}\left[\frac{d^{\frac{np}{q}}(ax^\alpha)}{dx^{\frac{np}{q}}}\right] = \frac{d^{\frac{p}{q}}}{dx^{\frac{p}{q}}}[a\frac{\Gamma(\alpha+1)}{\Gamma(\alpha-\frac{np}{q}+1)}x^{\alpha-\frac{np}{q}}]$$

$$= a\frac{\Gamma(\alpha+1)}{\Gamma(\alpha-\frac{np}{q}+1)}\frac{\Gamma(\alpha-\frac{np}{q}+1)}{\Gamma(\alpha-\frac{np}{q}-\frac{p}{q}+1)}x^{\alpha-\frac{np}{q}}$$

$$= a\frac{d^{\frac{(n+1)p}{q}}x^\alpha}{dx^{\frac{(n+1)p}{q}}} = \frac{d^{\frac{(n+1)p}{q}}(ax^\alpha)}{dx^{\frac{(n+1)p}{q}}}$$

**Theorem 4.4**

If the order $s$ is negative and the absolute value of $s$ is such $|s| \geq 1$, the fractional derivative as defined by the relation (4.1.2) is not unique and depends on $|[s]|$ arbitrary constants $c_k$.

Proof
If $s < 0$ and $|s| \geq 1 \Rightarrow k - [s] + 1 > 0$ in general. then the second term of the relation (4.1.2) is not equal to zero and the expression of a derivative $F(x)$ of order $s$ of $f(x) = ax^\alpha$ is

$$F(x) = \frac{d^s(ax^\alpha)}{dx^s} = a\frac{d^s x^\alpha}{dx^s}$$

$$= a\frac{\Gamma(\alpha+1)}{\Gamma(\alpha-s+1)}x^{\alpha-s} + \sum_{k=-|[s]|}^{-1} c_k \frac{x^{k-s}}{\Gamma(k-[s]+1)} \quad (4.1.11)$$

As the values of the constants $c_k$ are arbitrary, any other function $G(x)$ such

$$G(x) = F(x) + \sum_{k=-|[s]|}^{-1} c_k \frac{x^{k-s}}{\Gamma(k-[s]+1)} \quad (4.1.12)$$

is a derivative of order $s$ of $f(x)$ too.

### Theorem 4.5

If $s = -n \in \mathbb{Z}_-$, the derivative of order $s$ as defined by the relation (4.1.2) is the indefinite integral of order $n$.

Proof

If we replace $s = -n$ in the relation (4.1.2), we obtain

$$\frac{d^{-n}(ax^\alpha)}{dx^{-n}} = a\frac{d^{-n}x^\alpha}{dx^{-n}}$$

$$= a\frac{\Gamma(\alpha+1)}{\Gamma(\alpha+n+1)}x^{\alpha+n} + \sum_{k=n}^{-1} c_k \frac{x^{k-s}}{\Gamma(k+n+1)} \quad (4.1.13)$$

The last member of the relation (4.1.13) is equal to the result obtained by integrating $n$ times. The constants $c_k$ may be identified to the integration constants.

### Theorem 4.5

For $s = \frac{p}{q} \in \mathbb{Q}_-$, the expression of the derivative of order $s$, according to (4.1.2) is

$$\frac{d^{\frac{p}{q}}(ax^\alpha)}{dx^{\frac{p}{q}}} = a\frac{d^{\frac{p}{q}}x^\alpha}{dx^{\frac{p}{q}}}$$

$$= a\frac{\Gamma(\alpha+1)}{\Gamma(\alpha-\frac{p}{q}+1)}x^{\alpha+\frac{p}{q}} + \sum_{k=-[\frac{p}{q}]}^{-1} c_k \frac{x^{k-\frac{p}{q}}}{\Gamma(k-[\frac{p}{q}]+1)} \quad (4.1.14)$$

Proof

This relation is obtained directly from (4.1.2) by replacing $s = \frac{p}{q}$. We recall, according to the theorem 4.4, this fractional derivative is not unique for $\left|\frac{p}{q}\right| > 1$

Example

Applying the relation (4.1.14), we obtain

$$\frac{d^{-\frac{1}{4}}x^5}{dx^{-\frac{1}{4}}} = \frac{\Gamma(5+1)}{\Gamma\left(5+\frac{1}{4}+1\right)}x^{5-\frac{1}{4}} + \sum_{k=-\left\|\left[\frac{-1}{4}\right]\right\|}^{-1} c_k \frac{x^{k+\frac{1}{4}}}{\Gamma\left(k+\left[\frac{1}{4}\right]+1\right)}$$

$$= \frac{5!}{\Gamma(\frac{25}{4})}x^{\frac{21}{4}}$$

$$\frac{d^{-\frac{5}{3}}x^3}{dx^{-\frac{5}{3}}} = \frac{\Gamma(3+1)}{\Gamma\left(3+\frac{5}{3}+1\right)}x^{3+\frac{5}{3}} + \sum_{k=-1}^{-1} c_k \frac{x^{k+\frac{5}{3}}}{\Gamma\left(k+\left[\frac{5}{3}\right]+1\right)}$$

$$= \frac{3!}{\Gamma(\frac{18}{3})}x^{\frac{14}{3}} + c_{-1}$$

$$\frac{d^{-\frac{17}{7}}x^5}{dx^{-\frac{17}{7}}} = a\frac{\Gamma(5+1)}{\Gamma\left(5+\frac{17}{7}+1\right)}x^{5+\frac{17}{7}} + \sum_{k=-2}^{-1} c_k \frac{x^{k+\frac{17}{7}}}{\Gamma\left(k+\left[\frac{17}{7}\right]+1\right)}$$

$$= \frac{5!}{\Gamma(\frac{38}{7})}x^{\frac{31}{7}} + c_{-1}x + c_{-2}$$

Remark

If $s = \frac{p}{q} \in \mathbb{Q}_-$, it may be shown that we do not have exactly the relation

$$\underbrace{\frac{d^{\frac{p}{q}}}{dx^{\frac{p}{q}}} \cdot \frac{d^{\frac{p}{q}}}{dx^{\frac{p}{q}}} \cdots \frac{d^{\frac{p}{q}}}{dx^{\frac{p}{q}}}}_{n\ times}(ax^\alpha) = \frac{d^{\frac{np}{q}}(ax^\alpha)}{dx^{\frac{np}{q}}} \quad (4.1.15)$$

but it may be proved that we have an equality of the form

$$\frac{d^{\frac{np}{q}}(ax^\alpha)}{dx^{\frac{np}{q}}} = \underbrace{\frac{d^{\frac{p}{q}}}{dx^{\frac{p}{q}}} \cdot \frac{d^{\frac{p}{q}}}{dx^{\frac{p}{q}}} \cdots \frac{d^{\frac{p}{q}}}{dx^{\frac{p}{q}}}}_{n\ times}(ax^\alpha) + P(x) \quad (4.1.16)$$

in which $P(x)$ is a linear combination of power functions with arbitrary coefficients. As example, let us take $p = -1$ then we can show by mathematical induction that

$$\underbrace{\frac{d^{-\frac{1}{q}}}{dx^{-\frac{1}{q}}} \cdot \frac{d^{-\frac{1}{q}}}{dx^{-\frac{1}{q}}} \cdots \frac{d^{-\frac{1}{q}}}{dx^{-\frac{1}{q}}}}_{n\ times}(ax^\alpha) = a\frac{x^{\alpha+1}}{\alpha+1} \quad (4.1.17)$$

But we have

$$\frac{d^{-1}(ax^\alpha)}{dx^{-1}} = a\frac{x^{\alpha+1}}{\alpha+1} + c_{-1} \quad (4.1.18)$$

So

$$\frac{d^{-1}(ax^\alpha)}{dx^{-1}} = \underbrace{\frac{d^{-\frac{1}{q}}}{dx^{-\frac{1}{q}}} \cdot \frac{d^{-\frac{1}{q}}}{dx^{-\frac{1}{q}}} \cdots \frac{d^{-\frac{1}{q}}}{dx^{-\frac{1}{q}}}}_{n\ times}(ax^\alpha) + P(x)$$

with $P(x) = c_{-1}$ in our example. As the relation (4.1.17) differs from the relation (4.1.18) by the constant $c_{-1}$ we have a same property of the "semi-equality" introduced in paragraph 3.

### 4.2. Properties

#### 4.2.1. Linearity and Semi-linearity

### Theorem 4.5

If the order $s$ is a positive rational number ($s \in \mathbb{Q}_+$) or a negative rational number ($s \in \mathbb{Q}_-$) with $|s| < 1$, the fractional derivative of order $s$ (as defined in 4.1.2) has the property of linearity

$$\frac{d^s(ax^\alpha + bx^\beta)}{dx^s} = a\frac{d^s x^\alpha}{dx^s} + b\frac{d^s x^\beta}{dx^s} \quad (4.2.1)$$

Proof

If $s \in \mathbb{Q}_+$ or $s \in \mathbb{Q}_-$ with $|s| < 1$, we have

$$\sum_{k=-|[s]|}^{-1} c_k \frac{x^{k-s}}{\Gamma(k-[s]+1)} = 0$$

So

$$\frac{d^s a x^\alpha}{dx^s} = a \frac{d^s x^\alpha}{dx^s} = a \frac{\Gamma(\alpha+1)}{\Gamma(\alpha-s+1)} x^{\alpha-s}$$

Then

$$a \frac{d^s x^\alpha}{dx^s} + b \frac{d^s x^\beta}{dx^s}$$

$$= a \frac{\Gamma(\alpha+1)}{\Gamma(\alpha-s+1)} x^{\alpha-s} + b \frac{\Gamma(\alpha+1)}{\Gamma(\alpha-s+1)} x^{\beta-s}$$

$$= \frac{d^s(ax^\alpha + bx^\beta)}{dx^s}$$

**Theorem 4.8**

If the order $s$ is a negative rational number ($s \in \mathbb{Q}_-$) with $|s| > 1$, the fractional derivative, as defined in (4.1.2), doesn't have the property of linearity. In fact we have in this case

$$\frac{d^s(ax^\alpha + bx^\beta)}{dx^s}$$

$$= a \frac{d^s x^\alpha}{dx^s} + b \frac{d^s x^\beta}{dx^s} + \sum_{k=-|[s]|}^{-1} c_k \frac{x^{k-s}}{\Gamma(k-[s]+1)} \quad (4.2.2)$$

Proof

We have for $s \in \mathbb{Q}_-$ and $|s| > 1$ ($[s] \neq 0$)

$$\frac{d^s(ax^\alpha)}{dx^s} = a \frac{d^s x^\alpha}{dx^s}$$

$$= a \frac{\Gamma(\alpha+1)}{\Gamma(\alpha-s+1)} x^{\alpha-s} + \sum_{k=-|[s]|}^{-1} c_k \frac{x^{k-s}}{\Gamma(k-[s]+1)}$$

Then

$$a \frac{d^s x^\alpha}{dx^s} + b \frac{d^s x^\beta}{dx^s}$$

$$= a \frac{\Gamma(\alpha+1)}{\Gamma(\alpha-s+1)} x^{\alpha-s} + \sum_{k=-|[s]|}^{-1} c_{k1} \frac{x^{k-s}}{\Gamma(k-[s]+1)}$$

$$+ b \frac{\Gamma(\alpha+1)}{\Gamma(\alpha-s+1)} x^{\beta-s} + \sum_{k=-|[s]|}^{-1} c_{k2} \frac{x^{k-s}}{\Gamma(k-[s]+1)} \quad (4.2.3)$$

As the values of the constants $c_{k1}$ and $c_{k2}$ are arbitrary, we may identify $c_{k1} + c_{k2} = c_k$ then the relation (4.2.2) is deduced easily from (4.2.3) one. According to the definition in the paragraph 3, the relation (4.2.2) is a "semi-linearity" one.

### 4.2.2. Commutativity and Semi-commutativity

**Theorem 4.9**

For $(s_1, s_2) \in \mathbb{Q} \times \mathbb{Q}$ with $|s_1| < 1, |s_2| < 1$, the fractional derivative as defined by the relation (4.1.2) is commutative.

$$\frac{d^{s_1}}{dx^{s_1}} \frac{d^{s_2}(ax^\alpha)}{dx^{s_2}} = \frac{d^{s_2}}{dx^{s_2}} \frac{d^{s_1}(ax^\alpha)}{dx^{s_1}}$$

$$= \frac{d^{s_1+s_2}(ax^\alpha)}{dx^{s_1+s_2}} \quad (4.2.4)$$

Proof
For $|s| < 1$ we have $[s] = 0$ so

$$\frac{d^s ax^\alpha}{dx^s} = a \frac{d^s x^\alpha}{dx^s} = a \frac{\Gamma(\alpha+1)}{\Gamma(\alpha-s+1)} x^{\alpha-s}$$

Then

$$\frac{d^{s_1}}{dx^{s_1}} \frac{d^{s_2}(ax^\alpha)}{dx^{s_2}} = a \frac{\Gamma(\alpha+1)}{\Gamma(\alpha-s_2-s_1+1)} x^{\alpha-s}$$

$$= \frac{d^{s_1+s_2}(ax^\alpha)}{dx^{s_1+s_2}} = \frac{d^{s_1}}{dx^{s_1}} \frac{d^{s_2}(ax^\alpha)}{dx^{s_2}}$$

**Theorem 4.10**

If $(s_1, s_2) \in \mathbb{Q} \times \mathbb{Q}$ with $|s_1| \geq 1, |s_2| \geq 1$, the fractional derivative as defined by the relation (4.1.2) doesn't have the commutativity propriety.

Proof
From the definition (4.1.2), we obtain

$$\frac{d^{s_2}}{dx^{s_2}} \frac{d^{s_1}(ax^\alpha)}{dx^{s_1}} = a \frac{\Gamma(\alpha+1)}{\Gamma(\alpha-(s_1+s_2)+1)} x^{\alpha-(s_1+s_2)}$$

$$+ \sum_{k=-|[s_1]|}^{-1} c_k^{21} \frac{x^{k-(s_1+s_2)}}{\Gamma(k-[s_1]-s_2+1)}$$

$$+ \sum_{k=-|[s_2]|}^{-1} c_k^2 \frac{x^{k-s_2}}{\Gamma(k-[s_2]+1)}$$

$$\frac{d^{s_1}}{dx^{s_1}} \frac{d^{s_2}(ax^\alpha)}{dx^{s_2}} = a \frac{\Gamma(\alpha+1)}{\Gamma(\alpha-(s_1+s_2)+1)} x^{\alpha-(s_1+s_2)}$$

$$+ \sum_{k=-|[s_2]|}^{-1} c_k^{12} \frac{x^{k-(s_1+s_2)}}{\Gamma(k-[s_2]-s_1+1)}$$

$$+ \sum_{k=-|[s_1]|}^{-1} c_k^1 \frac{x^{k-s_1}}{\Gamma(k-[s_1]+1)}$$

So in general we have

$$\frac{d^{s_2}}{dx^{s_2}} \frac{d^{s_1}(ax^\alpha)}{dx^{s_1}} \neq \frac{d^{s_1}}{dx^{s_1}} \frac{d^{s_2}(ax^\alpha)}{dx^{s_2}} \quad (4.2.5)$$

On the equal footing as in paragraph 3, as the values of the constants $c_k$ are arbitrary, we may say that we have the semi-commutativity.

## 5. Conclusion and Comparison of the Two Definitions

The approach that we have introduced to study the fractional derivatives opens a new way and new insights concerning the definition and properties of the derivation operation, particularly in the possibility of an unified definition for both orders : positive and negative.

The two definitions that we have considered have their own advantages and disadvantages. The first one corresponding to the relation (3.1.2) seems to be more general but the existence of the infinite summation introduces some complexity which implies the study of the convergence of this series. The second one corresponding to the relation (4.1.2) seems to be simpler because it eliminates the existence of the infinite summation but it introduces some ambiguity particularly in the case of fractional negative orders as it can be seen in the relations (4.1.16).

The two definitions are exactly the same if the order is integer numbers: positive or negative.

We have introduced the notions of semi-equality, semi-linearity and semi-commutativity giving a broad sense to the relation (3.1.8). In this framework, the two definitions may be said to be "equivalent modulo the second term of this relation (3.1.8) ".

In our study, we have considered fractional order $s \in \mathbb{Q}$ and power functions defined on $\mathbb{R}_+$ and their linear combinations

In references [12],[13], we study another method to define real and complex order derivatives and integrals by using operator approach in the case of more general functions than power functions.